\documentclass{ifacmtg}
\usepackage[dcucite]{harvard}
\usepackage{amssymb,amsmath}
\usepackage{graphicx}

\begin{document}
\runauthor{S. B. Yuste}
\begin{frontmatter}
\title{Weighted average finite difference methods for fractional diffusion equations}
\author[Santos]{Santos B. Yuste\thanksref{Someone}}
\address[Santos]{Departamento de F\'{\i}sica, Universidad  de  Extremadura,
E-06071 Badajoz, Spain}
\thanks[Someone]{Partially supported by the Ministerio de
Ciencia y Tecnolog\'{\i}a (Spain) through Grant No. FIS2004-01399
and by the European Community's Human Potential Programme under
contract HPRN-CT-2002-00307, DYGLAGEMEM.}

\begin{abstract}
Weighted averaged finite difference methods for solving fractional diffusion
equations are discussed and  different formulae of the discretization of the
Riemann-Liouville derivative are considered. The stability analysis of the
different numerical schemes is carried out by means of a procedure close to the
well-known von Neumann method of ordinary diffusion equations. The stability
bounds are easily found and checked in some representative examples.
\end{abstract}

\end{frontmatter}

\section{Introduction}
The number of scientific and engineering problems involving fractional calculus
is very large and growing.  The applications range from control theory to
transport problems in fractal structures, from relaxation phenomena in
disordered media to anomalous reaction kinetics of subdiffusive reagents.
Recently,  a fractional Fokker-Planck equation has been proposed to describe
subdiffusive anomalous transport in the presence of an external field
\cite{MetBarKlaPRL,BarMetKlaPRE,MetzlerRev}. For the force-free case, the
equation becomes the fractional partial differential equation\cite{Balakrishnan,ScalasPRE,MetzlerRev,SchneiderWyss}
\begin{equation}
\label{subdeq} \frac{\partial }{\partial t} u(x,t)= K_\gamma
~_{0}D_{t}^{1-\gamma } \frac{\partial^2}{\partial x^2} u(x,t)
\end{equation}
where  $~_{0}D_{t}^{1-\gamma } $ is the fractional derivative defined by the
Riemann-Liouville operator, $K_\gamma$ is the diffusion coefficient and
$\gamma\in(0,1)$ is the anomalous diffusion exponent. There are many analytical
techniques for dealing with these fractional equations. But, as also with
ordinary (non-fractional) partial differential equations (PDEs), in many cases
the initial and/or boundary conditions and/or the external force are such that
the only reasonable option is to resort to numerical methods. However, these
methods are not as well studied as their non-fractional counterparts.

In this communication, some numerical methods for solving fractional partial
differential equations, which are very close to the well-known weighted average
(WA) methods of ordinary (non-fractional) PDEs, are considered. It will be
shown that the stability of the fractional numerical schemes can be analyzed
very easily and efficiently with a method close to von Neumman's (or Fourier's)
method for non-fractional PDEs.

\section{Fractional discretization formulae}
Two main steps will be considered to build numerical difference schemes for
solving fractional PDE's. In the first step one discretizes the ordinary
differential operators $\partial /\partial t$, $\partial^2/\partial x^2$ as
usual \cite{NumericalRecipes,MortonMayers}. This will be done in Sec.
\ref{sec:fracschemes}.  In the second step, one discretizes the
Riemann-Liouville operator:
\begin{equation}
\label{GLdefpura} ~_{0}D_{t}^{1-\gamma
}f(t)=\frac{1}{h^{(1-\gamma)}} \sum_{k=0}^{[t/h]}
\omega_k^{(1-\gamma)} f(t-kh) + O(h^p),
\end{equation}
where $[t/h]$ means the integer part of $t/h$. This formula is not
unique because there are many different valid choices for
$\omega_k^{(\alpha)}$ \cite{LubichCoef}.  Let $\omega(z,\alpha)$
be the generating function of the coefficients
$\omega_k^{(\alpha)}$, i.e.,
\begin{equation}
\label{} \omega(z,\alpha)=\sum_{k=0}^{\infty} \omega_k^{(\alpha)}
z^k .
\end{equation}
If the generating function is
\begin{equation}
\omega(z,\alpha)=(1-z)^\alpha
\end{equation}
then we get the backward difference formula of order $p=1$ (BDF1). This is also
called the backward Euler formula of order 1 or, simply the Gr\"unwald-Letnikov
formula. These coefficients are $\omega_k^{(\alpha)}=(-1)^k \binom{\alpha}{k} $
and can be evaluated recursively:
\begin{equation}
\omega_0^{(\alpha)}=1, \qquad
\omega_k^{(\alpha)}=\left(1-\frac{\alpha+1}{k}
\right)\omega_{k-1}^{(\alpha)}.
 \label{coefO1}
\end{equation}
The generating function for the backward difference formula of
order $p=2$ (BDF2) is
\begin{equation}
\label{}
   \omega(z,\alpha)=\left(\frac{3}{2}-2z+\frac{1}{2}
   z^2\right)^\alpha,
\end{equation}
and
\begin{equation}\label{}
  \omega(z,\alpha)=\left(\frac{11}{6}-3 z+\frac{3}{2}z^2-\frac{1}{3}z^3\right)^\alpha
\end{equation}
is the generating function for the backward difference formula of order $p=3$
(BDF3). Generating functions for higher-order BDF formulae can be found in
\cite{LubichCoef,Podlubny}. Another type of discretization formula is that of
Newton-Gregory of order $p$ (NG$p$) \cite{LubichCoef} whose generating function
is
\begin{equation}
\begin{aligned}\label{}
  \omega(z,\alpha)=&
(1-z)^\alpha \times \\
&\left[\Omega_0+\Omega_1(1-z)+\cdots \Omega_{p-1}\right]
\end{aligned}
\end{equation}
where the coefficients $\Omega_n$ are defined by
\begin{equation}\label{}
\sum_{n=0}^{\infty} \Omega_n (1-\xi)^n =\left(\frac{\ln \xi
}{\xi-1}\right)^\alpha.
\end{equation}

\section{Fractional difference schemes}
\label{sec:fracschemes} The notation $x_j=j \Delta x$, $t_m=m\Delta t$ and
$u(x_j,t_m)\equiv u_j^{(m)}\simeq U_j^{(m)}$ will be used with $U_j^{(m)}$
being the numerical estimate of $u(x,t)$ at the point $(x_j,t_m)$. In the
non-fractional weighted average method, the diffusion equation is replaced by:
\begin{equation}
\begin{aligned}
\label{ec:difufw} u_j^{(m+1)}=&u_j^{(m)}+ \lambda S
\left[{u_{j-1}^{(m)}-2u_j^{(m)}+u_{j+1}^{(m)}}\right]+\\
&(1-\lambda) S \left[{u_{j-1}^{(m+1)}-2u_j^{(m+1)}+u_{j+1}^{(m+1)}}\right]+\\
&T(x,t),
\end{aligned}
\end{equation}
where $\lambda$ is the weight factor,  $T(x,t)$  the truncation error
\cite{MortonMayers}, and  $S= D \Delta t/(\Delta x)^2$. Similarly, the
fractional equation is replaced by
\begin{equation}
\label{ec:met1}
\begin{aligned}
 u_j^{(m+1)}=&u_j^{(m)}+ \lambda S
\left[{ \overset{\bullet} u_{j-1}^{(m)}-2\overset{\bullet}u_j^{(m)}+\overset{\bullet}u_{j+1}^{(m)}}\right]+\\
&(1-\lambda) S \left[{\overset{\bullet}u_{j-1}^{(m+1)}-2\overset{\bullet}u_j^{(m+1)}+\overset{\bullet}u_{j+1}^{(m+1)}}\right]+\\
&T(x,t),
\end{aligned}
\end{equation}
where $\overset{\bullet}u_j^{(m)}\equiv ~_{0}D_{t}^{1-\gamma
}u(x_j,t_m)$.
 Inserting Eq.\ \eqref{GLdefpura}  into Eq.\ (\ref{ec:met1}), neglecting the truncation
error, and rearranging the terms, we finally get the fractional WA
difference scheme
\begin{equation}
\label{erme}
\begin{aligned}
 U_j^{(m+1)}=& U_j^{(m)}+ (1-\lambda) S \sum_{k=0}^{m+1}
\omega_k^{(1-\gamma)}
\left[U_{j-1}^{(m+1-k)}- \right. \\
& \left. 2U_j^{(m+1-k)}+U_{j+1}^{(m+1-k)}\right] + \lambda \,S \times \\
& \sum_{k=0}^{m} \omega_k^{(1-\gamma)}
\left[U_{j-1}^{(m-k)}-2U_j^{(m-k)}+U_{j+1}^{(m-k)}\right]\; ,
\end{aligned}
\end{equation}
where
\begin{equation}\label{sgammas}
S= \frac{K_\gamma {\Delta t}}{{h^{1-\gamma}(\Delta x)^2}} \; .
\end{equation}
For $\lambda=1$ we recover the fractional explicit method discussed in
\cite{YusteAcedoFracExpli,YusteAcedoFracExpli2}. The method is implicit for $\lambda\neq 1$. For
$\lambda=0$ one gets the (fractional) fully implicit method, and for
$\lambda=1/2$  the (fractional) Crank-Nicholson method. Because the estimates
$U_j^{(m)}$ of $u(x_j,t_m)$ are made at the times $m \Delta t$, $m=1,2,\ldots$,
and because the evaluation of $~_{0}D_{t}^{1-\gamma} u(x_j,t)$ by means of
\eqref{GLdefpura} requires knowing $u(x_j,t)$  at the times $n h$, $n=0,
1,2,\ldots$, it is natural to choose $h=\Delta t$. In this case,
\begin{equation}
\label{sgam} S= K_\gamma \frac{\Delta t^\gamma}{(\Delta x)^2}\; .
\end{equation}

Before tackling Eq.~\eqref{erme} seriously, one must first know under which
conditions, if any, the integration algorithm is stable.

\section{Stability analysis}
\begin{figure}
\begin{center}
\includegraphics[width=0.95 \columnwidth]{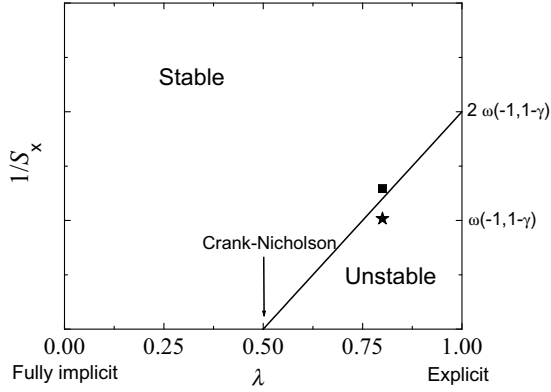}
\end{center}
\caption{Stability phase diagram for the weighted average methods versus the
weight factor $\lambda$. The line corresponds to Eq.~\eqref{1StimesMain}. The
square corresponds to the case shown in Fig.~\ref{fig5} and the star
corresponds to the cases shown in Figs.~\ref{fig6} and \ref{fig7}.
\label{fig1}}
\end{figure}

\begin{figure}
\begin{center}
\includegraphics[width=0.95 \columnwidth]{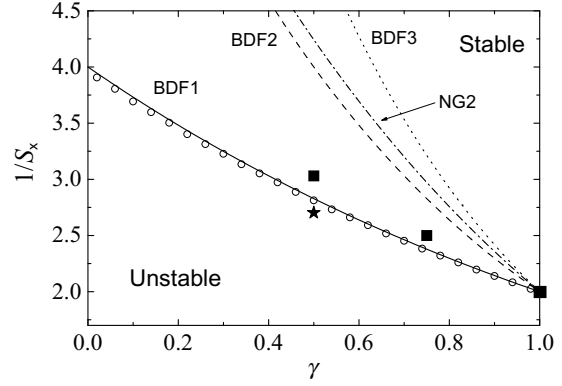}
\end{center}
\caption{Stability bound for the explicit method  versus the anomalous exponent
$\gamma$. The lines BDF1, BDF2, BDF3 and NG2 correspond to the theoretical
stability bounds of the explicit method when the BDF1, BDF2, BDF3 and NG2
discretization formulae of the Riemann-Liouville derivative are used. The
circles are numerical estimates of the BDF1 stability bound
\protect\cite{YusteAcedoFracExpli,YusteAcedoFracExpli2}. The squares correspond to the cases shown
in Fig.~\ref{fig3} and the star is the case of Fig.~\ref{fig4}. \label{fig2}}
\end{figure}

The stability analysis of the integration difference scheme \eqref{erme}  will
be carried out by means of the method used in \cite{YusteAcedoFracExpli,YusteAcedoFracExpli2}.
Following the von Neumann ideas, one studies  the stability of a single
\emph{subdiffusive} mode of the form $U_j^{(m)}=\zeta_m e^{i q j \Delta x}$.
This mode will  be stable as long as  the $\zeta_m$ stay bounded for $m\to
\infty$. Proceeding as \cite{YusteAcedoFracExpli,YusteAcedoFracExpli2} one finds that  a fractional
WA method is stable as long as $1/S\ge 1/S_\times$ where
\begin{equation}\label{1StimesMain}
1/S_\times= 2(2\lambda-1) \omega(-1,1-\gamma).
\end{equation}
This is the main result of the present work. From \eqref{1StimesMain} one sees
that a WA method is stable for any value of $S$ if $\lambda\le 1/2$ because the
generating function $\omega(z,1-\gamma)$ for $z=-1$ is positive ($S$ is always
positive: see Eq.\ \eqref{sgammas}). For $\lambda>1/2$, $S_\times$ is positive,
and there exist values of $1/S$ smaller than $1/S_\times$, so that the
fractional WA methods are unstable for these cases.  Figure~\ref{fig1} shows
the stability phase diagram for the WA methods. For $\lambda=1$ one recovers
the stability bound $1/S_\times= 2\omega(-1,1-\gamma)$ for the fractional
explicit methods discussed in \cite{YusteAcedoFracExpli,YusteAcedoFracExpli2}. The stability bounds
of these explicit methods versus the anomalous diffusion exponent $\gamma$ for
several Riemman-Liouville derivative discretization formulae are shown in
Fig.~\ref{fig2}. Note that the stability region shrinks when the order $p$ of
the discretization formula increases.

The fully implicit numerical method ($\lambda=0$) and the Crank-Nicholson
method ($\lambda=1/2$) have been considered, together with their stability, by
some authors \cite{Serna,Lopez,LubichST,McLean1,McLean2}. However, their
analysis is far more complex than that presented here. Also, to the best of my
knowledge, neither generic implicit WA methods (arbitrary $\lambda$) nor the
explicit method ($\lambda=1)$ have been  considered previously.  In this vein,
it is interesting to note that an \emph{explicit} method for solving the
fractional diffusion equation written in the Caputo form has recently been
proposed by Ciesielski and Leszczynski \cite{Ciesielski}.

\section{Check of the stability bound}
\begin{figure}
\begin{center}
\includegraphics[width=0.95 \columnwidth]{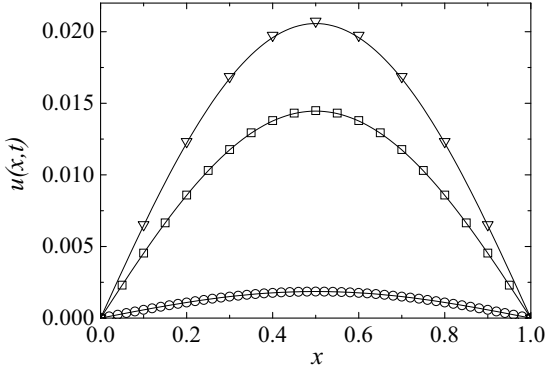}
\end{center}
\caption{Solution of the subdiffusion equation \eqref{subdeq} with
absorbing boundary conditions $u(0,t)=u(1,t)=0$ and initial
condition $u(x,0)=x(1-x)$. The symbols correspond to the
BDF1-explicit numerical solution and the lines correspond to the
exact analytical solution. The solution is shown for the time
$t=0.5$ for $\gamma=0.5$, $S=0.33$, $\Delta x=1/10$ (triangles), $\gamma=0.75$, $S=0.4$, $\Delta x=1/20$  (squares) and $\gamma=1$, $S=0.5$, $\Delta x=1/50$ (circles). These cases are marked by squares in Fig.
\protect\ref{fig2} \label{fig3}}
\end{figure}

\begin{figure}
\begin{center}
\includegraphics[width=0.95 \columnwidth]{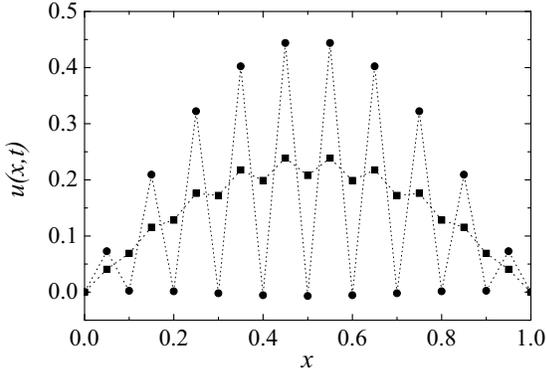}
\end{center}
\caption{BDF1-explicit numerical solution $u(x,t)$ for the same problem as in
Fig.~\ref{fig3} but with $\gamma=1/2$, $K_\gamma=1$, $S=0.37$, and $\Delta
x=1/20$ after 150 time steps (squares) and 200 time steps (circles). The lines
are plotted as a visual aide. This case corresponds to the point marked by a
star in Fig. \protect\ref{fig2}. \label{fig4}}
\end{figure}

The stability of the fractional difference schemes of
Sec.~\ref{sec:fracschemes} will be checked here by applying them
to solve Eq.~\eqref{subdeq} with the initial condition
$u(x,t=0)=x(1-x)$ and the absorbing boundary conditions
$u(0,t)=u(1,t)=0$. The exact analytical solution of Eq.\
(\ref{subdeq}) is easily found by the method of separation of
variables:
\begin{equation}
\begin{aligned}
\label{uexact} u(x,t)=&\frac{8}{\pi^3} \sum_{n=0}^\infty
\frac{1}{(2n+1)^3} \sin[(2n+1)\pi x] \times \\
& E_\gamma[-K_\gamma (2n+1)^2\pi^2 t^\gamma]\; .
\end{aligned}
\end{equation}
where $E_\gamma$ is the Mittag-Leffler function \cite{MetzlerRev,MainardiGorenfloJCAM}. In Fig.\
\ref{fig3} we compare this exact solution with the results provided by the BDF1
\emph{explicit} integration scheme for anomalous diffusion exponents
$\gamma=0.5$, $\gamma=0.75$, and $\gamma=1$ for time $t=0.5$ and $K_\gamma=1$.
The values of $\Delta x$ used were $\Delta x=1/10$, $1/20$, and $1/50$ with
$S=0.33$, $0.4$, and $0.5$, respectively. These values of $S$ are marked by
squares in Fig.~\ref{fig2}. They are inside the stable region, which is
confirmed by the well-behaved numerical solutions shown in Fig.~\ref{fig3}.
However, in Fig.\ \ref{fig4} the BDF1-explicit numerical solution for $S=0.37$
and $\gamma=1/2$ shows a characteristic unstable behaviour. This is the
expected behavior because $1/0.37$ is smaller than $1/S_\times=2
\omega(-1,1-\gamma)=2^{2-\gamma}$ for $\gamma=1/2$. This case is marked by a
star in Fig.~\ref{fig2}.

\begin{figure}
\begin{center}
\includegraphics[width=0.95 \columnwidth]{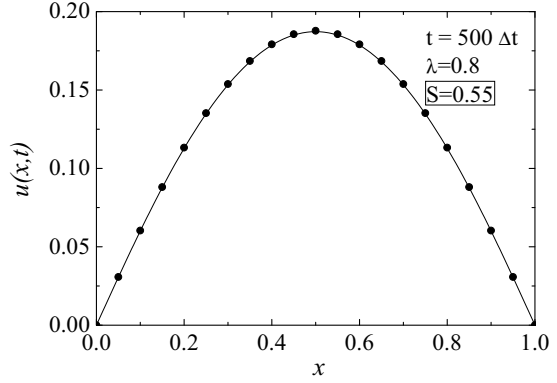}
\end{center}
\caption{Numerical solution $u(x,t)$ obtained by means of the weighted average
method with weight factor $\lambda=0.8$  for the same problem as in
Fig.~\ref{fig3} but with $\gamma=1/2$, $K_\gamma=1$, $S=0.55$, and $\Delta
x=1/20$ after 500 time steps (circles). The line is the exact analytical
result. This case corresponds to the point marked by a square in Fig.
\protect\ref{fig1}. \label{fig5}}
\end{figure}

\begin{figure}
\begin{center}
\includegraphics[width=0.95 \columnwidth]{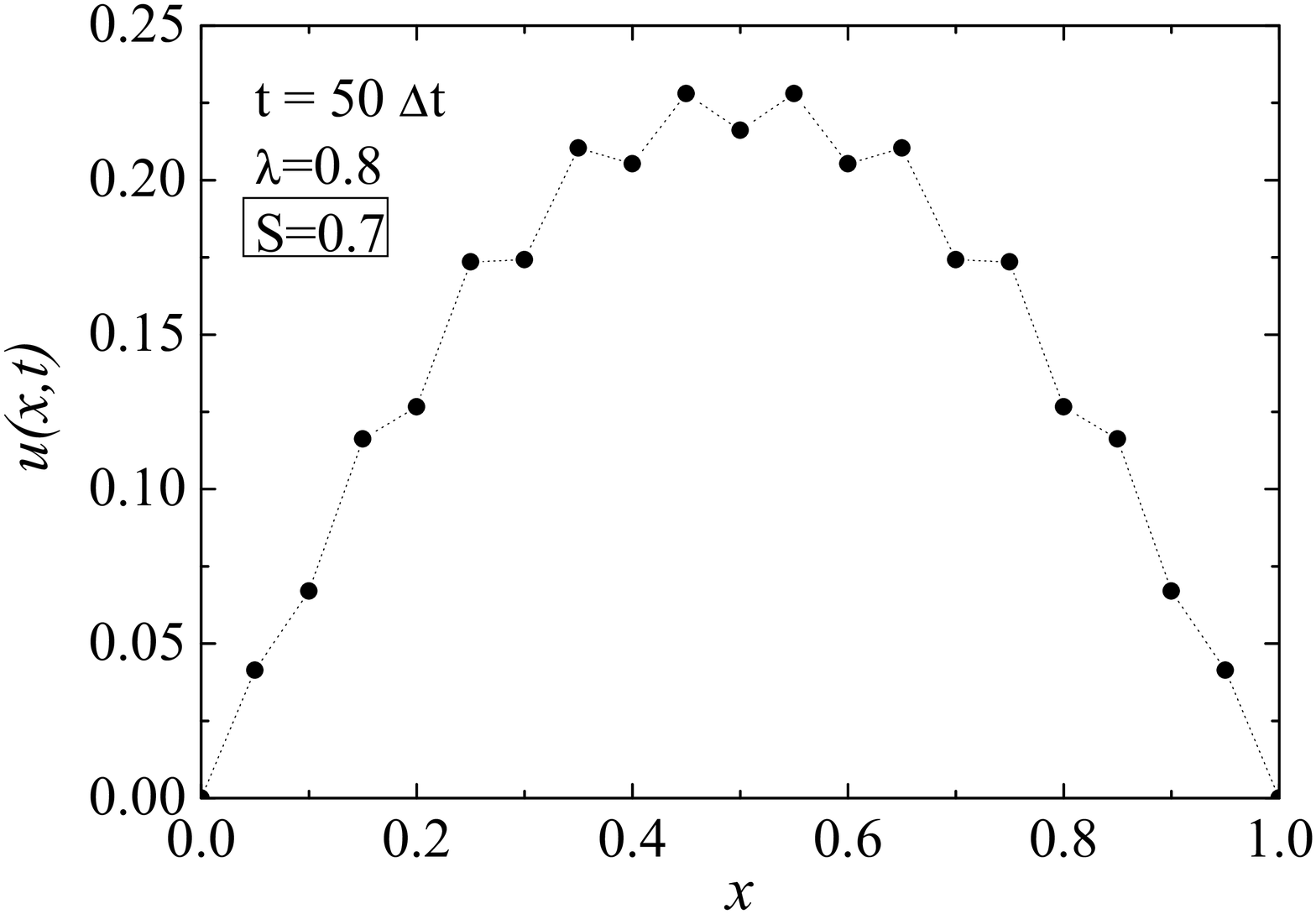}
\end{center}
\caption{Numerical solution $u(x,t)$ for the same problem as in Fig.~\ref{fig3}
obtained by means of the weighted average method with weight factor
$\lambda=0.8$  for $\gamma=1/2$, $K_\gamma=1$, $S=0.7$, and  $\Delta x=1/20$
after 50 time steps. The lines are plotted as a visual aide. This case
corresponds to the point marked by a star in Fig. \protect\ref{fig1}.
\label{fig6}}
\end{figure}

\begin{figure}
\begin{center}
\includegraphics[width=0.95 \columnwidth]{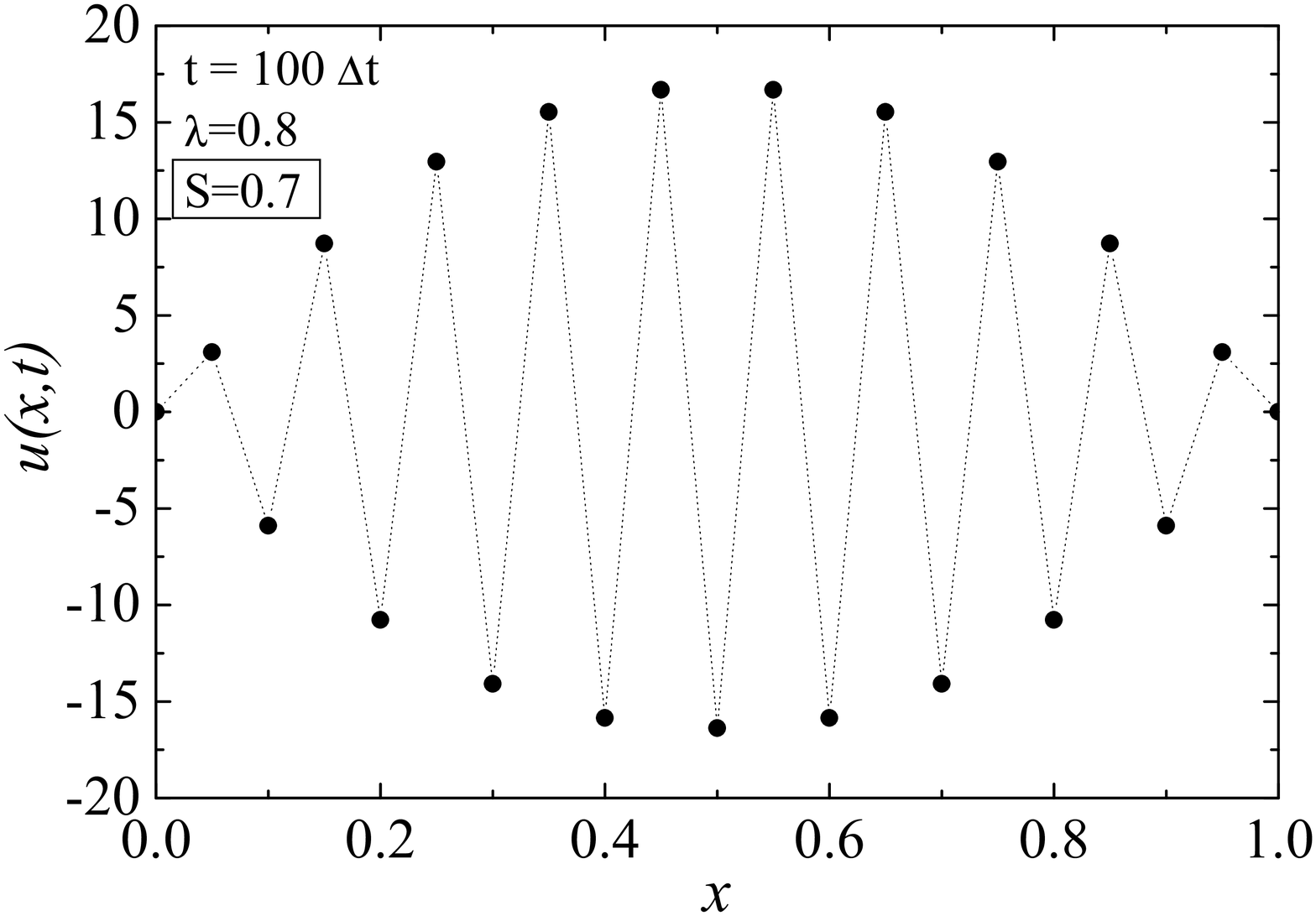}
\end{center}
\caption{Numerical solution $u(x,t)$ for the same problem as in Fig.~\ref{fig3}
obtained by means of the weighted average method with weight factor
$\lambda=0.8$ for $\gamma=1/2$, $K_\gamma=1$, $S=0.7$, and $\Delta x=1/20$
after 100 time steps. The lines are plotted as a visual aide. This case
corresponds to the point marked by a star in Fig. \protect\ref{fig1}.
\label{fig7}}
\end{figure}

Figures \ref{fig5}, \ref{fig6}, and \ref{fig7} show the numerical integration
results for the WA implicit method with $\lambda=0.8$, $\gamma=1/2$ and two
values of $S$. For $S=0.55$ one has $1/S>1/S_\times=1.2\times 2^{1/2}$ (this
case corresponds to the square in Fig.~\ref{fig1}) and the WA method must be
stable. This is confirmed in Fig.~\ref{fig5}. However,
$1/S<1/S_\times=1.2\times 2^{1/2}$ for $S=0.7$ so that the WA method must be
unstable in this case, which is confirmed in Figs.~\ref{fig6} and \ref{fig7}
(this case is marked by a star in Fig.~\ref{fig1}).


\section{Final remarks}
The truncation error $T(x,t)$ of Eq.~\ref{ec:met1} can be estimated in the same
way as for non-fractional equations \cite{MortonMayers}. One gets:
 \begin{equation}
\label{txt14}
\begin{aligned}
T_j^m=&O(h^p)+O(\Delta t)^2+O(\Delta x)^2+
\left(\frac{1}{2}-\lambda\right) O(\Delta t)+\\
& (1-\lambda)O\left(\frac{\omega_{m+1}^{(1-\gamma)}}{h^{1-\gamma}}
\right) \frac{\partial^2}{\partial x^2}u(x_j,t_{1/2})
\end{aligned}
\end{equation}
where $T_j^m$ is the truncation error at $x_j$ and $t_m+\Delta t/2$:
$T_j^m\equiv T(x_j,t_{m+1/2})$. From this expression some conclusions may be
drawn:

\noindent $\bullet$  If $h=\Delta t$, it is useless to employ discretization
formulae for the Riemann-Liouville derivative of order $p$ higher than two
because of the unavoidable presence of an $O(\Delta t)^2$ term.

\noindent $\bullet$ The undiserable low-order term proportional to $\Delta t$
is present for all WA methods with $\lambda\neq 1/2$. Therefore,  the chosen
value of $\lambda$ (as long as $\lambda\neq 1/2$) does not improve the
precision of the WA method (although it matters for the stability of the
numerical scheme). The value $\lambda=1/2$ is special because it removes the
$O(\Delta t)$ term. It leads to the (fractional) Crank-Nicholson method.

\noindent $\bullet$ The last term does not appear for non-fractional
discretization schemes: it is characteristic of fractional methods. It becomes
negligible for large $m$. In fact, for $m$ large enough, the quantity
$\omega_{m+1}^{(1-\gamma)}/h^{1-\gamma}$ becomes of order of, or smaller than,
$O(\Delta t)^2$. The particular value of $m$ for which this happens depends on
the discretization formula of the Riemann-Liouville derivative. However, for
the first integration steps ($m$ small) this term, and consequently, the
truncation error, is large unless the initial curvature of $u(x,t)$,
$\partial^2 u(x,t)/\partial x^2$, is small. These difficulties disappear for
$\lambda=1$, i.e., they are absent in the explicit method.  This suggests a
practical integration procedure in which the first integration steps are
performed by means of the explicit method, and the subsequent steps are carried
out by means of, say, the fractional Crank-Nicholson method.

Of the numerical integration schemes for fractional PDE's considered here, the
Crank-Nicholson method appears to be the most promising  because it is always
stable and is second order accurate in $\Delta t$ and $\Delta x$ (provided that
the discretization formula for the Riemann-Liouville derivative is second order
accurate in $h=\Delta t$).

It has been shown that the stability of the WA methods can be studied by means
of a von-Neumann-type stability analysis that is surprisingly simple and
accurate. The fractional integration schemes and stability analysis discussed
here can be easily extended to $d$-dimensional fractional diffusive equations
and fractional wave equations.


\end{document}